\documentstyle[12pt,twoside]{article}
\input xy
\xyoption{all}

\font\capi=cmbx10  scaled\magstep 2     
\font\para=cmsy10   scaled\magstep 3    
\font\petita=cmcsc10 scaled\magstep 1   

\font\es=msbm10 scaled\magstep 1
\font\ess=msbm10

\def\Rational{\mbox{\es Q}}
\def\Complex{\mbox{\es C}}
\def\Proj{\mbox{\es P}}
\def\proj{\mbox{\ess P}}

\begin{document}

\baselineskip=16pt

\def\leaderfill{\leaders\hbox to 1em {\hss.\hss}\hfill}

\centerline {\LARGE{\bf Numerical Bounds of Canonical Varieties}}

$$
\begin{array}{l}
\mbox{\rm {\bf Miguel A. BARJA}}\\
\end{array}
$$

\vglue1.7truecm

\noindent {\capi {\para \S} 0. Introduction}

\vglue.3truecm

Let $X$ be a minimal, complex, projective, Gorenstein variety of dimension $n$. We
say that $X$ is $canonical$ if for some (any) desingularization
$\sigma : Y \longrightarrow X$, the map associated to the
canonical linear series $\vert K_Y \vert$ is birational.

We note $K_{X}$ for the canonical divisor of $X$ and
$\omega_{X}={\mathcal {O}}_X (K_X)$ the canonical
sheaf. Let $p_g=h^0(X,\omega_X)$, $q=h^1(X,{\mathcal {O}}_X)
$. There are several known bounds for $K_X^n$ depending on $p_g$, the most
general one being the bound $K_X^n \geq (n+1)p_g+d_n$ ($d_n$ constant)
given by Harris (\cite{H}). Bounds including other invariants are
known for canonical surfaces, $K_S^2 \geq 3p_g+q-7$ (\cite{J},
\cite{De}), and for surfaces and threefolds fibred over curves (\cite{O},
\cite{X2}).

In this paper we prove some results for canonical surfaces and threefolds.
In the case of canonical surfaces there are some known results which
show that under some additional hypotheses, the bound $K_S^2 \geq
3 p_g+q-7$ can be considerably improved (see Remark 2.2). We give here some other special
cases (Remark 2.2) for which is not sharp and prove (Theorem 2.1) that,
in fact, $K_S^2 = 3 p_g +q-7$  only if $q=0$ whenever
$p_g(S)\geq 6$.

Canonical surfaces with $K_S^2=3p_g-7$ are known to exist and
classified (\cite{AK}).Then we can hope that a good bound for canonical
surfaces including the irregularity should be of type $K_S^2 \geq
3 p_g +aq-7$, $a >1$. Since for $q=1$ it is known (\cite{Ko2}) that
$K^2_S \geq 3 p_g$, $a$ should be 7, although unfortunately
examples of low $K_S^2$ (with $q \geq 2$) are not known.

In the case of canonical threefolds we prove that $K_X^3 \geq 4p_g+6q-32$.
In particular, we prove that the results of Ohno for canonical fibred
threefolds are not sharp.

\hrule

\baselineskip=12pt

\vglue.5truecm

\noindent
{\footnotesize Partially supported by CICYT PS93-0790 and HCM project n.
ERBCHRXCT-940557}

\newpage

\baselineskip=16pt

We use basically a result on quadrics containing irreducible varieties
due to Reid (\cite{R1}) and several techniques originated in \cite{X1}
and developed by Konno (\cite{Ko2}, \cite{Ko3}, \cite{Ko4}) for the
study of the slope of fibred surfaces. In particular we include in an
Appendix the dimension 3 version of the relative hyperquadrics method
used by Konno in \cite{Ko4}.

After this manuscript was written, the author was informed that
theorem 2.1 was known yet to K. Konno (unpublished).

The author wants to thank his advisor, professor Juan C. Naranjo,
for fruitful conversations and continuous support.

\vglue.7truecm

\noindent {\capi {\para \S} 1. A general inequality}

\vglue.3truecm

We need the following result due to Reid (\cite{R1}, p. 195).

\bigskip

\noindent {\bf Lemma 1.1.} \
{\it
Let $\Sigma \subseteq {\Proj}^N$ be an irreducible
variety spanning ${\Proj}^N$ of dimension $w$. Then
$$h^0 {\mathcal {J}}_{\Sigma ,{\proj}^N}(2) \leq {{N-w+2}\choose{2}}-
\mbox{\rm min} \{ \mbox{\rm deg} \Sigma,2(N-w)+1\}.$$
}

Then we have an immediate consequence.

\bigskip

\noindent {\bf Proposition 1.2.} \
{\it
Let $X$ be a normal projective variety of general type and dimension $n$.
Let $L\in Div(X)$,
${\mathcal {L}}={\mathcal {O}}_X(L)\in Pic X$ and $\varphi$ the rational map
associated to $\mathcal {L}$. Assume $\varphi$ is birational; then

{\rm (a)} $h^0(X,{\cal O}_X(2L)) \geq (n+2)[h^0(X,{\mathcal {O}}_X(L))-\frac{n+1}{2}]$

{\rm (b)} If equality holds in {\rm (a)} then

\qquad {\rm (i)} $\Sigma := \varphi(X)$ is contained in a minimal degree
variety of ${\Proj}^{h^0(X,{\mathcal {L}})-1}$ of dimension $n+1$ obtained as the
intersection of quadrics containing $\Sigma$.

\qquad {\rm (ii)} $\Sigma \subseteq {\Proj}^{h^0(X,{\cal L})-1}$ is linearly
and quadratically normal.

\qquad {\rm (iii)} $Bs \vert L \vert = Bs \vert 2L \vert$.

\qquad {\rm (iv)} If $Bs \vert L \vert = \emptyset$ and $p,q \in X$
then $\vert L \vert$ separates $p$ and $q$ if and only if so does
$\vert 2L \vert$.
}

\medskip

\noindent {\petita Proof:}

We can always consider
$$\xymatrix{
{\bar X} \ar[d]_{\sigma} \ar@{>>}[dr]^{\bar {\varphi}} \\
X \ar@{-->}[r]_-{\varphi} & \Sigma \subseteq {\Proj}^r
}$$
where $r=h^0({\mathcal {L}})-1$, $\bar X$ is smooth, $\sigma$ is birational and
$\bar{\varphi}$ is defined by the moving part $M$ of the linear system
$\vert {\sigma}^{\ast} (L) \vert$, which has no base point.

By construction we have ${\bar{\varphi}}^{\ast}
{\mathcal {O}}_{{\proj}^r}(1)=
{\mathcal {O}}_{\bar X}(M)$ and $2M \leq $ moving part of $\vert {\sigma}^{\ast}
(2L) \vert$. Then, since $X$ is normal and $\sigma$ has connected fibres
$$\begin{array}{rl}
h^0(X,{\mathcal {O}}_X(2L))&=h^0(X,{\sigma}_{\ast}{\sigma}^{\ast}
{\mathcal {O}}_X(2L))=\\
&=h^0(\bar X,{\sigma}^{\ast}{\mathcal {O}}_X(2L))\geq h^0(\bar X,
{\mathcal {O}}_{\bar X}(2M))=\\
    &=h^0(\bar X, {\bar \varphi}^{\ast} {\mathcal {O}}_{{\proj}^r}(2))=\\
        &=h^0(\Sigma,{\bar \varphi}_{\ast} {\bar \varphi}^{\ast}
           {\mathcal {O}}_{{\proj}^r}(2)) \geq h^0(\Sigma, {\mathcal {O}}_{\Sigma}(2))
\end{array} \qquad \ \, (1) $$

Now if we consider
$$0 \longrightarrow H^0 {\mathcal {J}}_{\Sigma,{\proj}^r}(2) \longrightarrow
H^0{\mathcal {O}}_{{\proj}^r}(2) \buildrel f \over \longrightarrow
H^0 {\mathcal {O}}_{\Sigma} (2)$$
Lemma 1.1 gives
$$\begin{array}{rl}
h^0(\Sigma,{\mathcal {O}}_{\Sigma}(2)) &\geq \mbox {\rm dim Im}f \geq
{{r+2}\choose{2}}-{{r+2-n}\choose{2}}+
\mbox{\rm min} \{ \mbox{\rm deg} \Sigma,2(r-n)+1\}\\
&=(n+2)[r-\frac{n-1}{2}]=(n+2)[h^0({\cal O}_X(L))-\frac{n+1}{2}]
\end{array} \ (2)$$
if deg$\Sigma \geq 2(r-n)+1$. If $H_i\ (i=1,\dots,n)$ are general
hyperplanes in ${\Proj}^r$ and ${\Sigma}_k=\Sigma \cap H_1 \cap \dots
\cap H_{n-k}$ is a general section of $\Sigma$ of dimension $k$ we have
that ${\Sigma}_2$ is an irreducible surface of general type and then
(\cite{Be1}, p. 115):
$$\mbox{\rm deg} \Sigma = \mbox{\rm deg}{\Sigma}_2 \geq 2(r-n+2)-1>
2(r-n)+1.$$

This proves (a).

Assume from now on that equality holds in (a). In particular equality
must hold at every step of (1) and (2). Then $f$ is an epimorphism
and $h^1 {\mathcal {J}}_{\Sigma,{\proj}^r}(2)=0$. Since $h^1 {\mathcal {J}}_{\Sigma,
{\proj}^r}(1)$ is always zero we have (ii). Moreover we have
$$\begin{array}{ccccc}
S^2H^0{\mathcal {O}}_{{\proj}^r}(1)&\longrightarrow &
S^2H^0 {\mathcal {O}}_{\Sigma}
(1)\cong& S^2H^0(\bar X, {\sigma}^{\ast}{\mathcal {L}})
\cong &S^2H^0(X,{\mathcal {L}})\\
\downarrow & & & &\downarrow {\alpha}  \\
H^0{\mathcal {O}}_{{\proj}^r} (2) &\buildrel f \over \longrightarrow
&H^0 {\mathcal {O}}_{\Sigma}(2)
\cong &H^0 (\bar X, {\sigma}^{\ast}{\mathcal {O}}_X (2L))
\cong &H^0(X,{\mathcal {O}}_X(2L))
\end{array}$$
and hence $\alpha$ is an epimorphism and (iii) follows immediately.

In order to prove (iv), consider local trivializations of $\mathcal {L}$ at
$p$ and $q$. For $\alpha,\beta \in H^0({\cal L})$ we confuse $\alpha,\beta$ with
their local expressions at these trivializations.

We need

\bigskip

\noindent {\bf Claim.} If $Bs \vert L \vert = \emptyset$ then
$\vert L \vert $ does not separate $p$ and $q$ if and only if for all
$\alpha,\beta \in H^0(L)$, $\left \vert \begin{array}{cc}
\alpha(p)&\beta(p)\\ \alpha(q)&\beta(q) \end{array} \right \vert=0$.

\medskip

\noindent {\petita Proof} of the Claim:

Let $\beta \in H^0({\cal L})$ be such that
$\beta(p)=0$. Since $p$ is not a base point of $\vert L \vert$ there
exists $\alpha \in H^0({\cal L})$ such that $\alpha(p) \not= 0$. Then,
from $\left \vert \begin{array}{cc}
\alpha(p)&0\\ \alpha(q)&\beta(q) \end{array} \right \vert=0$
we get $\beta(q)=0$ and then $\beta$ does not separate $p$ and $q$.

Assume there exist $\alpha,\beta \in H^0({\cal L})$ such that
$\alpha (p)=a$,
$\alpha (q)=b$, $\beta (p)=\bar a$, $\beta (q)=\bar b$ and
$a \bar b - b \bar a \not= 0$. Let $\sigma=a \beta - \bar a \alpha \in
H^0({\cal L})$. Then clearly $\sigma$ separates $p$ and $q$.
$\quad \Box$

\medskip

If $\vert 2L \vert$ does not separate $p$ and $q$ then trivially so does
not $\vert L \vert$.

Assume $\vert L \vert$ does not separate $p$ and $q$. Since
$S^2H^0({\cal L})
\longrightarrow H^0({\cal L}^{\otimes 2})$ is surjective
for every $\alpha, \beta \in H^0(2L)$,
$\alpha=\sum a_{ij}s_is_j$, $\beta=\sum b_{ij}s_is_j$, $s_i \in
H^0({\cal L})$.
Since $\vert L \vert$ has no base point and does not separate
$p$ and $q$ we can take $\bar s \in H^0({\cal L})$ such that $\bar s(p)=a \not=
0$, $\bar s(q) =b \not= 0$. Since by the claim we have
$\left \vert \begin{array}{cc} s_i(p)&a\\s_i(q)&b \end{array} \right \vert
= 0$ for every $s_i$ we can define $\lambda_i=\frac{s_i(p)}{a}=
\frac{s_i(q)}{b}$. Then $\alpha(p)=\sum a_{ij} \lambda_i \lambda_j
a^2$, $\alpha(q)=\sum a_{ij}\lambda_i \lambda_j b^2$, $\beta (p)=\sum
b_{ij}\lambda_i \lambda_j a^2$, $\beta(q)=\sum b_{ij} \lambda_i
\lambda_j b^2$, and then
$$\left \vert \begin{array}{cc} \alpha(p)&\beta(p)\\\alpha(q)&
\beta(q) \end{array} \right \vert =a^2 b^2 \left \vert
\begin{array}{cc} \sum a_{ij} \lambda_i \lambda_j & \sum b_{ij}
\lambda_i \lambda_j\\ \sum a_{ij}\lambda_i \lambda_j & \sum
b_{ij} \lambda_i \lambda_j \end{array} \right \vert =0 $$
and hence, by the claim, $\vert 2L \vert$ does not separate $p$ and $q$.

For the proof of (i) we refer again to \cite{R1} p. 195. If we call
$\Sigma_0=\Sigma \cap H_1 \cap \dots \cap H_n$ we have that $\Sigma_0$
is a set of $d=\mbox{\rm deg} \Sigma \geq 2(r-n)+3$ points in
${\Proj}^{r-n}$. Proof of Lemma 1.1 (cf. \cite{R1} p.195)
shows that if we consider
$$H^0{\mathcal {O}}_{{\proj}^r}(2) \buildrel f \over \longrightarrow
H^0 {\mathcal {O}}_{\Sigma}(2)$$
$$H^0 {\mathcal {O}}_{{\proj}^{r-n}} (2) \buildrel f_0 \over \longrightarrow
H^0 {\mathcal {O}}_{\Sigma_0}(2)$$
then dim Im$f \geq {{r+2}\choose {2}} - {{r+2-n}\choose {2}} +$
dim Im$f_0 \geq {{r+2}\choose{2}}
-{{r+2-n}\choose{2} } +$ min $\{d,2(r-n)+1\}$. Under our
hypothesis equality holds and then we have that $\Sigma_0$ is a set of
$d$ points in ${{\Proj}^{r-n}}$ imposing exactly $2(r-n)+1$ conditions
on quadrics. Then $\Sigma_0$ is contained in a rational normal curve
$\Gamma$ intersection of the quadrics containing $\Sigma_0$. Let $T_k$
be the intersection of quadrics of ${\Proj}^{r-n+k}$ containing
$\Sigma_k$. We have $T_k \subseteq T_{k+1}\cap H_{n-k}$ and hence
$\Gamma = T_0 \subseteq T_n \cap H_1 \cap \dots \cap H_n$. Then $T_n$ has
an irreducible component $W$ containing $\Sigma$ of dimension at least
$n+1$. But then
$$h^0 {\mathcal J}_{W,{\proj}^r}(2) = h^0 {\mathcal {J}}_{\Sigma,
{\proj}^r}(2) =
{{r-n}\choose {2}} $$
since $\Sigma \subseteq W \subseteq T_n$. Again applying Lemma 1.1
to $W$, if $w=$ dim$W\geq n+2$:
$$h^0 {\mathcal {J}}_{W,{\proj}^r}(2)\leq {{r-n}\choose{2}} -1.$$

So dim$W=n+1$ and, since $W \cap H_1 \cap \dots \cap H_n=\Gamma$, $W$ is
a variety of minimal degree in ${\Proj}^r$. Since such varieties are
always intersection of quadrics we have in particular $W=T_n$.
$\quad \Box$

\vglue.7truecm

\noindent {\capi {\para \S} 2. Canonical surfaces}

\vglue.3truecm

As a consequence of Proposition 1.1 we get the following result for
minimal canonical surfaces.
The first part is a well known fact (cf. \cite{De}, \cite{J}).

\bigskip

\noindent {\bf Theorem 2.1.} \
{\it
Let $S$ be a minimal canonical surface. Then

{\rm (a)} $K_S^2 \geq 3 p_g+q-7$.

{\rm (b)} Assume $p_g(S)\geq 6$.
If $K_S^2=3p_g+q-7$, then $q=0$.
}

\medskip

\noindent {\petita Proof:}

(a) By the pluri-genus formula, we have $h^0(2K_S)=K_S^2+\chi
({\cal O}_S)$.
Hence the inequality $K_S^2 \geq 3p_g+q-7$ follows immediatly from
Proposition 1.2 (see \cite{De}, \cite{J}).

(b) In order to prove the statement we need first some properties of
surfaces lying on the border line; let $\Sigma=\varphi (S) \subseteq
{\Proj}^{p_g-1}$.

\bigskip

\noindent {\bf Claim 1.} If $K_S^2=3p_g+q-7$ then

(i) $\Sigma$ lies in a threefold $Z$ of minimal degree.

(ii) $\vert K_S \vert$ is base point free.

(iii) $\vert K_S \vert$ does not separate $p,q \in S$ (possibly
infinitely near) if and only if so does not $\vert 2 K_S \vert$.

(iv) $q=0$ or $q \geq 3$.

(v) If dim$\,Sing \Sigma=1$ and $K_S^2 \geq 13$ then the one dimensional
components of $Sing \Sigma$ are double lines.

\medskip

\noindent {\petita Proof} of Claim 1:

(i), (ii) and (iii) are direct consequence of Proposition
1.1 and the fact that $\vert 2 K_S \vert$ has no base points if $p_g \geq 4$
(\cite{Bo}).

(iv) If $q \not=0$ and $q \leq 2$, $K_S=3p_g+q-7$, then
$K_S^2 < 3 \chi {\cal O}_S$ and the canonical map of $S$ can not be
birational (cf. \cite{Ho}).

(v) Assume dim$\,Sing \Sigma =1$. Let $D$ be a one dimensional component
of $Sing \Sigma$. The canonical map $\varphi$ is not an embedding over
$D$. Since $K_S^2 \geq 10$ and since, by (iii) points which are not
separated by $\vert K_S \vert $ are those which are not separated
by $\vert 2 K_S \vert$ we can apply Reider's Theorem (see \cite{Re}).
Let $q \in D$ be a general point of $D$ and let
$p_1,p_2 \in S$ (possibly infinitely near) such that $\varphi(p_1)=
\varphi (p_2)=q$. By Reider's Theorem  we have that there exists an
effective divisor $E$ passing through $p_1,p_2$ and verifying
$$(K_SE,E^2)=(0,-2),(1,-1),(2,0).$$

Since $q \in D$ is general, $E$ can not be contracted by $\varphi$.
Since $\varphi(p_1)=\varphi(p_2)$, we have that $\mbox{deg}(\varphi
_{\vert E}) \geq 2$ and hence that $(K_SE,E^2)=(2,0)$ holds.

Note that moving $q \in D$ the curve $E$ can not move because then we
would have the surface $S$ covered by curves of genus at most two and this
is impossible since $S$ is canonical. So we must have $\varphi(E)=D$
(set-theoretically), deg$\varphi_{\vert E}=2$ and that $D$ is a line
${\Proj}^{p_g-1}$. In order to prove that it is a double line, assume
that, for $q \in D$ general, there are 3 points $p_1,p_2,p_3$ (possibly
infinitely near) such that $\varphi (p_i)=q$. Consider the 0-cycle
$Z=p_1+p_2+p_3$. Then $Z$ is clearly in special position with respect
to $\vert 2K_S \vert$ (following Reider's notation in \cite{Re}). Then
if $K_S^2 \geq 4 \mbox{deg}Z +1=13$, we have a decomposition $K_S=M+F$
as in \cite{Re}, where $F$ is an effective curve passing through $Z$
(and so that $K_SF \geq 3$) with the following numerical conditions:

$$M^2 \geq F^2, MF \geq 0, M^2>0$$
$$F^2 \leq \frac{(MF)^2}{(M^2)}$$
$$K_S^2=M^2+2MF+F^2 \geq 13$$
$$MF \, \mbox{is even}.$$

An easy computation show that these conditions have no solution.
Then a general point $q \in D$ has two pre-images through $\varphi$.
$\quad \Box$

\medskip

It is a well known fact that the only possibilities for a threefold $Z$
of minimal degree in ${\Proj}^{p_g -1}$ are

(A) $Z={\Proj}^3\ (p_g=4)$.

(B) $Z$ is a cone over the Veronese surface ($p_g=7$).

(C) $Z$ is a smooth quadric in ${\Proj}^4 \ (p_g=5)$.

(D) $Z$ is a scroll of type ${\Proj}_{a,b,c}$, \
$0 \leq a \leq b \leq c$,\
$2 \leq a+b+c=p_g-3$.

\bigskip

\noindent {\bf Claim 2.} If $p_g(S) \geq 6$ and $K_S^2=3p_g+q-7$,
then $q=0$.

\medskip

\noindent {\petita Proof} of Claim 2:

We will often use the following results on the irregularity of a surface
of general type:

\begin{itemize}

\item[(*)] (Beauville, \cite{Be3}): $p_g \geq 2q-4$. If $S$ is canonical,
then $p_g(S) \geq 2q-3$.

\item[(**)] (Xiao, \cite{X2}): if $f:S \longrightarrow {\Proj}^{1}$ is a
linear pencil of general fibre $F$, then $q \leq \frac{1}{2}(g(F)+1)$.

\item[(***)] (Konno, \cite{Ko5}): with the same notation as in (**), is
$r$ is the rank of the subsheaf of $f_{*} \omega_{S}$ generically
generated by its global sections, then $q \leq \mbox{min}\{g(F)-r, r+1\}$.

\end{itemize}

\medskip

Assume first that $p_g(S)=7$ and that $Z$ is a cone over the Veronese
surface. Choosing hyperplane sections of $Z$ passing through the vertex
and decomposing as two components as sections of the Veronese surface,
we have a decomposition:

$$K_S \sim 2D+D^{'}$$
where $\vert D \vert$ is a linear pencil on $S$, $DD^{'} \geq 1$
(by the connectedness of canonical divisors) and that
$K_SD^{'}=0$ ($D^{'}$ is contracted by $\varphi$ to the vertex
of $Z$; possibly $D^{'}=0$). Then, the only possibility is that
$(K^2_S,q)=(18,4)$ and hence that $g(D) \leq 7$. But using (***) we
obtain than $r \leq 3$ wich is impossible since $\varphi (D)$ is not
a plane curve.

\medskip

Assume $Z$ is a scroll. Consider
$$\xymatrix{
\bar S \ar[d]_\sigma \ar[r] & \bar \Sigma \subseteq \bar Z \ar[d]
\ar[r] & {\Proj}^1\\
S \ar[r] & \Sigma \subseteq Z
}$$
where $\bar Z$ is the desingularization of $Z$. Let $\bar \alpha:
\bar S \longrightarrow {\Proj}^1$ be the induced fibration and $\bar G$
be a general fibre. Note that, by construction
$(\varphi \circ \sigma)_{\vert \bar G}:\bar G \longrightarrow {\Proj}^{p_g-1}$
induces on $\bar G$ a base point free sublinear system of $\vert
K_{\bar G} \vert$ and that $(\varphi \circ \sigma )(\bar G) \subseteq
{\Proj}^2 \cong T$, where $T$ is a general ruling of $Z$.

Note that the singularities of $(\varphi \circ \sigma) (\overline G)$,
for $\bar G$ general, lie on
$Sing Z$ (produced by the base points of $\vert \sigma(\bar G)\vert$ on
$S$) or on $Sing \Sigma \cap T$. If $a+b+c \geq 2$ (we only exclude the
case
$Z={\Proj}^3$ which is Case (A)) then $p_g \geq 6$ and $K_S^2 \geq
11+q \geq 14$ if $q \not= 0$. Then, if $Sing \Sigma$ has one dimensional
components, they must be double lines by Claim 1. Moreover we can assume that
they are transversal to the general ruling. Since any such line in $Z$
corresponds to an epimorphism ${\mathcal {O}}_{{\proj}^1}(a)\oplus
{\mathcal {O}}_{{\proj}^1}(b) \oplus {\mathcal {O}}_{{\proj}^1}(c)
\longrightarrow {\mathcal {O}}_{{\proj}^1}(1)$, under the assumption
$a+b+c \geq 4$ ($p_g \geq 7$) we have
that the lines transversal to the ruling cut a general plane $T$ in
points which are on a line $\ell \subseteq T$. Then we can proceed as
follows.

Assume first $Z$ is smooth, i.e. $1 \leq a \leq b \leq c$. We have then
that $S=\bar S$, $\sigma(\bar G)=G$ and $\varphi(G)$ is a plane curve
of degree $d=2g(G)-2$ with only double points as singularities, lying
all of them on a line if $p_{g} \geq 7$. Let $m$ be the number of
such double points (possibly $m=0$). Then we have:

$$d \geq 2m,$$
$$d={\bar G}{\cal O}(1)=GK_S=GK_S +G^2=2g-2,$$
$$2g-2=\frac{1}{2}(d-1)(d-2)-m,$$
which forces $d \leq 5$ and hence $g=g(G) \leq 3$. By (**) we get a
contradiction.

If $p_g(S)=6$ and $Z$ is a smooth scroll of ${\Proj}^{5}$ we can proceed
as follows. Since ${\bar G} \subset {\Proj}^{2}$ we have that $r=3$ and
hence that $g(G) \geq 6$ by (***).

Take a general
hyperplane section of $\Sigma$ and $Z$. We get an irreducible curve
$\bar C$ lying on a smooth ruled surface $V$ of minimal degree in
${\Proj}^{4}$. Let $h$, $f$ be the hyperplane divisor class and the
fibre divisor class in $V$. We have that $h^2 = \mbox{\rm deg}V=p_g-3=3$
and that ${\bar C}=\alpha h +\beta f$ with $\alpha \geq 1$,
$\beta \geq -\alpha$ (obeserve that $V$ is ${\Proj}({\cal O}_{{\Proj}^1}
\oplus {\cal O}_{{\Proj}^1}(1))$ embeded in ${\Proj}^4$ by $C_0+2f$
($C_0^2=-e$), and hence $e=1$, ${\bar C}\sim \alpha C_0+(\beta+2 \alpha)f$
and so $\beta+2 \alpha \geq e \alpha=\alpha$,
following \cite{Be1} page 69 and \cite{Ha} Corollary V.2.18).
Let $C \in \vert K_S \vert$ be the smooth curve lying over $\bar C$. Using
that $K_V=-2h+(p_g-5)f$ we get
$$\begin{array}{rl}
d&=\mbox{deg}{\bar G}={\bar G}H=f{\bar C}= \alpha\\
d&=GK_S=K_SG+G^2=2g(G)-2 \geq 10\\
K_S^2&=\mbox{\rm deg} (\bar C)=\bar C h=\alpha (p_g-3)+\beta
 \geq 2\alpha \\
K_S^2&=3p_g+q-7=11+q
\end{array}$$
which is impossible.

\smallskip

Assume dim$\,SingZ=1$, i.e. $0=a=b<c$.
Take a general section $\Gamma$
of $\Sigma$ containing $SingZ$. $\Gamma$ corresponds to a section
$\vert K_S \vert \ni C=cG+L$ where $L$ is the component of the sublinear
system containing $SingZ$ (possibly $L=0$).

We have then, since $p_g \geq 2q-3$
$$\frac{7}{2}c+5 \geq 3p_g+q-7=K_S^2=cKG+KL\geq cKG.$$

Then, using $c \geq 3$, $cG^2\leq K_SG$ and evenness of $K_SG+G^2$
we get that, in any case $2p_a(G)-2=K_SG+G^2\leq 6$. Then
$g(\bar G)\leq p_a(G)\leq 4$. Again by (**) we get $q \leq 2$
and hence $q=0$.

\smallskip

Finally assume dim$\,Sing Z=0$, i. e. $0=a <b \leq c$.
We have a Hirzebruch surfce $\Sigma_{c-b} \subset {\Proj}^{p_g-2}$ as a
smooth hyperplane section of $Z$ embedded by $\vert C_0+cf \vert$
($C_0^2=-e=-(c-b)$). If we take now hyperplane sections of $Z$ passing
through its vertex and having exactly the above decomposition when
meeting
$\Sigma_{c-b}$,
then we produce a decomposition

$$\begin{array}{rl}
K_S&=cD+D'+D''\\
K_SD&\geq 1\\
K_SD'&\geq 1\\
K_SDïï&=0
\end{array}$$
such that $\vert D \vert$ is a linear pencil on $S$ with $\varphi (D)$
a plane curve. Note that if $b=c$ (then $e=0$), then $\vert D' \vert$
also moves on $S$. A case by case computation of the cases $c \geq 4$,
$c=3$, $c=2,b=1$, and $c=b=2$ shows that we always get a contradiction
by using (**) and (***).

\smallskip

We get then that
the only possibilities for $S$ with $q \not= 0$ occur when $p_g \leq 5$.
$\quad \Box$

\bigskip

\noindent {\petita Remark 2.2.} Part (b) of Theorem 2.1
shows that inequality
$K_S^2 \geq 3 p_g +q - 7$ is not sharp if $p_g>>0$.
Since surfaces with $K_S^2 =
3 p_g - 7$ are known to exist (and are completely understood, see
\cite{AK}), it seems that a sharp bound should look like $K_S^2
\geq 3 p_g +aq -7$, with $a > 1$. There are several partial results in
this direction:

\medskip

\noindent
(i) Let $alb:S \longrightarrow alb(S)$ be the Albanese map of $S$. As
a direct consequence of the study of the slope of fibrations, Konno
(\cite{Ko2}) shows that, if dim$\,alb(S)=1$ then $K_S^2 \geq 3 p_g
+7q -7$.

\smallskip

\noindent
(ii) In the same paper Konno proves that if the cotangent sheaf of $S$
is nef then $K_S^2 \geq 6 \chi {\mathcal {O}}_S =6 p_g-6q+6$ which is better
than $K_S^2 \geq 3p_g+q-7$ if $p_g >>q$.

\smallskip

\noindent
(iii) Note that even if dim$\,alb(S)=2$ but there exists a fibration
$\pi: S \longrightarrow B$ with $b=g(B)\geq 2$ we have $K_S^2 \geq
3 p_g +2q -7$. Indeed, for a general fibration we have $K_S^2 \geq
\lambda \chi {\mathcal {O}}_S + (8 - \lambda)(b-1)(g-1)$ ($g=g(F)$, $F$ smooth
fibre of $\pi$). Note that if $S$ is canonical $g \geq 3$. Under our
hypothesis $\pi \not= alb$ and then Xiao (\cite{X1}) proves that
$\lambda \geq 4$. Finally note that since $b + g \geq q$ (\cite{Be3} we
have
$$\begin{array}{rl}
(b-1)(g-1) &\geq (b-1)+(g-1) \geq q-2 \quad {\mbox {\rm if}}\ b \geq 3\\
{\mbox {\rm and}}\
(b-1)(g-1) &\geq (b-1)+(g-1)-1 \geq q-3 \quad {\mbox {\rm if}}\ b =2
\end{array}$$

But if $b=2$ and $(b-1)(g-1)=q-3$ we have $q=b+g$. Again by
\cite{Be3} we can say that $S=B \times F$ with $b=g(B)=2$. This is not
possible if $S$ is canonical. Finally we can apply (*)
and we get the desired bound.

\smallskip

\noindent
(iv) Let $C \in \vert K_S \vert$; then we have
$$0 \rightarrow H^0({\mathcal {O}}_S) \rightarrow H^0(\omega_S)
\rightarrow H^0(C,{\omega_S}_{\vert C}) \rightarrow
H^1({\mathcal {O}}_S) \buildrel \rho_C \over
\rightarrow H^1(\omega_S) \rightarrow \dots$$

Note that the above sequence is self-dual and then we can consider
$\rho_C \in Sym \,{\Complex}^q$. The correspondence $H^0(S,\omega_S)
\buildrel \Phi \over \longrightarrow Sym \,{\Complex}^q$ is clearly linear
since it is induced by the natural map $H^0(S,\omega_S) \otimes
H^1 ({\mathcal {O}}_S) \longrightarrow H^1(S,\omega_S)$. Then, if $p_g \geq
{{q+1}\choose{2}}$ there must exist $C \in \vert K_S \vert$
such that $\rho_C=0$. For such $C$ we have $h^0(C,{\omega_S}_{\vert C})=p_g+q-1$.

Assume $C$ to be irreducible. Since the linear system $\vert K_S
{\vert}_{\vert C}$ is birational we can apply ``Clifford plus" (\cite{R1}
p. 195) and get
$$p_g+q-1 = h^0(C,{\omega_S}_{\vert C}) \leq \frac{1}{3} (K_S^2+4)$$
and hence $K_S^2 \geq 3p_g + 3q -7$.

\vglue.7truecm

\noindent {\capi {\para \S} 3. Canonical threefolds}

\vglue.3truecm

\noindent

\noindent {\bf Theorem 3.1.}\
{\it
Let $T$ be a canonical threefold. Then
$$K_T^3 \geq 4p_g +6q -32$$
}

\medskip

\noindent {\petita Proof:}

Since $T$ is canonical (in particular, $T$ is minimal),
$K_T$ is nef and big and hence by
the general Kawamata-Viehweg Theorem (\cite{K2} Thm. 2.17) and
Proposition 1.2 we get
$$\frac{1}{2} K_T^3 - 3 \chi {\mathcal {O}}_T =
\chi_T (\omega_{T}^{\otimes 2})=
h^0(T,\omega_{T}^{\otimes 2}) \geq 5(h^0(T,\omega_T)-2)\kern.9truecm
(1)$$
and hence
$$K_T^3 \geq 4 p_g + 6(h^2 ({\mathcal {O}}_T) -
h^1 ({\mathcal {O}}_T)) -14 $$

Assume $h^2 ({\mathcal {O}}_T) \geq 2 h^1 ({\mathcal {O}}_T) -3$; then we get
$$K_T^3 \geq 4 p_g +6q -32  \qquad \qquad \qquad \qquad \qquad \qquad \qquad
\qquad \qquad \quad \qquad \qquad (2)$$
and then the Theorem is proved under this hypothesis.

From now on we assume $h^2({\mathcal {O}}_T) \leq 2 h^1({\mathcal {O}}_T) -4$; then
by \cite{Be1} Lemma X.7 and \cite{Be2} Proposition 1 we
obtain the existence of a fibration $\pi: T \longrightarrow B$ where $B$
is a smooth curve of genus $b \geq 2$.

Let $F$ be a general fibre of $\pi$. Since $K_T + F_{\vert F} =K_F$ we
have that the general fibre is a smooth canonical minimal surface (note
that $K_T$ is nef so in particular it is $\pi$-nef).

Then we can apply the results of Ohno (\cite{O}) and state that
(Main Theorem 2):
$$K_T^3-6(b-1)K_F^2=K_{T/B}^3 \geq 4 (\chi {\mathcal {O}}_B \chi
{\mathcal {O}}_F - \chi {\mathcal {O}}_T )\qquad
\qquad \qquad \qquad \quad(3)$$
except for a finite number of exceptions. We have
$$\begin{array}{rl}
K_T^3 &\geq 2 (b-1) [3 K_F^2-2 \chi {\mathcal {O}}_F ] + 4 p_g -4
(h^2 ({\mathcal {O}}_T) - h^1 ({\mathcal {O}}_T)) -4 \geq \\
&\geq 2 (b-1) [3 K_F^2-2 \chi {\mathcal {O}}_F ] + 4 p_g -4 q+12 \end{array}$$
since we are assuming $h^2 {\mathcal {O}}_T \leq 2 h^1({\mathcal {O}}_T) -4$. Note
that since $F$ is canonical
$$3 K_F^2 -2 \chi {\mathcal {O}}_F \geq 7 p_g (F) +5q(F) -23 \geq
5(q(F)-1)$$
and
$$2(b-1)[3K_F^2 - 2 \chi {\mathcal {O}}_F] \geq 10(b-1)(q(F)+1)\geq
10(q(F)+b)-10$$
since $b \geq 2$, $q(F)\geq 0$.

Note also that from the Albanese maps associated to $F \hookrightarrow
T \longrightarrow B$ we get $q(F) + b \geq q(T) = q$ and so
$$K_T^3 \geq 4 p_g +6q+2$$
which is stronger than we wanted.

Finally we must deal with the exceptions of Main Theorem 2 in
\cite{O}. Notice that since $F$ is a canonical surface we must have,
by Section 2, $K_F^2 \geq 3 p_g(F) + q(F) -7$.
From this, only a few exceptional cases hold. We divide
them in three cases (following \cite{Ko1} a canonical surface
verifying $p_g(F)=6$, $q(F)=0$, $K_F^2=3p_g(F)-6=12$ is classified
in two types according its canonical image is contained in a threefold of
$\Delta$-genus 0 or 1). In all of them we will prove $K_{T/B}^3
\geq 4(\chi {\mathcal {O}}_B \chi {\mathcal {O}}_F -\chi {\mathcal {O}}_T)
$. Then the same argument as above works.

\bigskip

\noindent {\sl Case 1.-} $p_g(F)=4,5$.

We use the results of the relative hyperquadrics method of the Appendix.
If ${\mathcal {E}} = \pi_{\ast} {\omega}_{T/B}$ and we consider the
relative canonical image of $T$:
$$\xymatrix{
T \ar@{-->}[r]^-\psi \ar[d] &Y \subseteq {\Proj}_B({\mathcal {E}})=Z
\ar[dl]^\varphi \\
B}$$

Then formula (A.2) gives
$$K_{T/B}^3 \geq (2 p_g(F)-4)
(\chi {\mathcal {O}}_B \chi {\mathcal {O}}_F- \chi {\mathcal {O}}_T)-
2 \mbox{deg}K-2 \ell (2)$$
where $K=\varphi_{\ast} {\mathcal {J}}_{Y,Z}(2)$. Note that since
$T$ is Gorenstein, $\ell (2)=0$ ([F]).

If $p_g(F)=4$, $K=0$ and
$$K_{T/B}^3 \geq 4(\chi {\cal O}_B \chi {\mathcal {O}}_F - \chi
{\mathcal  O}_T) $$
which produces, as in (3)
$$K_T^3 \geq 4 p_g+6q+2 \qquad  \qquad \qquad \qquad \qquad \qquad
\qquad \qquad \qquad \qquad \quad (4)$$

If $p_g(F)=5$ then rk$K=1$ and deg$K=x$ for some relative hyperquadric
$Q \equiv 2 L_{\mathcal {E}} -xF$ containing $Y$ (see proof of Lemma
A.4). Lemma A.5 of the Appendix
gives that deg$K=x \leq \frac{2}{3}$deg$\mathcal {E}$ since rk$Q \geq 3$.
Then from the proof of Corollary A.2 we get
$$\begin{array}{rl}
K_{T/B}^3 &\geq 2(p_g(F)+1)\mbox{\rm deg}{\mathcal {E}} - 6
(\chi {\mathcal {O}}_B \chi {\mathcal {O}}_F - \chi {\mathcal {O}}_T) -2 \mbox{\rm deg}
K\geq\\
&\geq \frac{32}{3} \mbox{\rm deg}{\mathcal {E}} -6(\chi {\mathcal {O}}_B \chi
{\mathcal {O}}_F - \chi {\mathcal {O}}_T)\geq \\
&\geq \frac{14}{3} (\chi {\mathcal {O}}_B \chi {\mathcal {O}}_F -
\chi {\mathcal {O}}_T)\geq \\
&\geq 4 (\chi {\mathcal {O}}_B \chi {\mathcal {O}}_F -
\chi {\mathcal {O}}_T)
\end{array}$$
which gives again (4).

\bigskip

\noindent {\sl Case 2.-} $p_g(F)=6,7$, $q(F)=0$ and $K_F^2=3p_g-7$ or
$p_g(F)=6$, $q(F)=0$, $K_F^2=3p_g-6$ and the canonical image of $F$ is
contained in a threefold of $\Delta$-genus 0, intersection of the quadrics
containing it.

\smallskip

Consider again the relative canonical image of $T$.

$$\xymatrix{
T \ar@{-->}[r]^-\psi \ar[d]_\pi &Y \subseteq
{\Proj}_B({\mathcal {E}})=:Z \ar[dl]^
\varphi \\
B}$$

If $A \in PicB$ is ample enough we have an epimorphism
$$H^0({\mathcal {J}}_{Y,Z}(2 L_{\mathcal {E}} \otimes \varphi^{\ast} (A))
\longrightarrow H^0({\mathcal {J}}_{F, {\proj}^{p_g-1}}(2))$$

Let $W$ be the horizontal irreducible component of the base locus
of the linear system given by the sections of $H^0({\mathcal {J}}_{Y,Z}
(2 L_{\mathcal {E}} \otimes {\varphi}^{\ast}(A)))$. Since under our
hypothesis intersections of quadrics containing $F$ is a threefold of
minimal degree (see \cite{AK} and \cite{Ko1}) $W$ is a
fourfold fibred over $B$ by threefolds of minimal degree. Let
$\widetilde W$ be a desingularization of $W$.

We want to relate the invariants of $\pi:T \longrightarrow B$ with those
of $\Phi : \widetilde W \longrightarrow B$. In \cite{Ko3}, Konno
gives a general method for this. We refer there for details. Let $H$ be
the pull-back of the tautological divisor of $Z$ to $\widetilde W$.

\bigskip

\noindent {\bf Lemma 3.2.}
{\it

{\rm (a)} $\Phi_{\ast} {\mathcal {O}}_{\widetilde W} (H)=\pi_{\ast}
\omega_{T/B}$.

{\rm (b)} $\mbox{\rm deg} \Phi_{\ast} {\mathcal {O}}_{\widetilde W}(2H)
=H^4+4 \mbox{\rm deg} \pi_{\ast} \omega_{T/B}$.

{\rm (c)} $K_{T/B}^3 \geq 2 H^4 + 2(\chi {\mathcal {O}}_B \chi
{\mathcal {O}}_F - \chi {\mathcal {O}}_T) $.
}

\medskip

\noindent {\petita Proof:}

(a) Follows directly from the construction of $\widetilde W$ and $H$.

(b) Note that the formula we want to prove is invariant under the
change of $H$ by $H + \Phi^{\ast}(A)$, $A \in PicB$. So we can
assume $\vert H \vert$ is base point free and hence get a smooth ladder
$\widetilde W=W_4 \supseteq W_3 \supseteq W_2 \supseteq W_1 \supseteq
W_0$ (i.e., $W_i$ is smooth and $W_i \in \vert H_{\vert W_{i+1}}
\vert $).
Notice that $W_2$ is a ruled surface over $B$. By induction one easily
proves that
$$ \forall i \geq 0 \ \forall m \geq 1 \ \forall n \geq 0 \quad
R^m \Phi_{\ast} {\mathcal {O}}_{W_i} (n H_{\vert W_i})=0$$
and hence that
$$\mbox{\rm deg} \Phi_{\ast} {\mathcal {O}}_{W_i} (2 H_{\vert W_i})=
\mbox{\rm deg} \Phi_{\ast} {\mathcal {O}}_{\widetilde W} (H)+
\mbox{\rm deg}\Phi_{\ast} {\mathcal {O}}_{W_{i-1}}(2 H_{i-1}).$$

Finally note that deg$\Phi_{\ast} {\mathcal {O}}_{W_0} (2H)=H^4$.

(c) The natural map $0 \longrightarrow \Phi_{\ast}
{\mathcal {O}}_{\widetilde W}(2H) \longrightarrow \pi_{\ast}
\omega_{T/B}^{\otimes 2}$ has a torsion
cokernel since it is an isomorphism at a general fibre. Then the result
follows calculating deg$\pi_{\ast} w_{T/B}^{\otimes 2}$ as in proof
of Corollary A.2 and applying (b).
$\quad \Box$

\medskip

In order to finish {\sl Case 2} note that, since part (c) of Lemma holds,
it is enough to prove that $H^4 \geq \mbox{\rm deg} \Phi_{\ast}
{\cal O}_{\buildrel \sim \over W}(H)$.

\bigskip

\noindent {\bf Claim:} Let $X$ be a smooth variety and $f:X\longrightarrow
B$ a filtration onto a smooth curve. Let $D \in Div(X)$ be a nef divisor
and let ${\cal E}=f_{\ast}{\cal O}_X(D)$. Then $D^n \geq \mbox{\rm deg}
f_{\ast} {\cal O}_X(D)$.

\noindent {\petita Proof} of the Claim: It follows easily by induction
from \cite{Ko3} Lemma 2.1.
$\quad \Box$

\bigskip

\noindent {\sl Case 3.-} $p_g(F)=6$, $q(F)=0$, $K_F^2=12$ and the
canonical image of $F$ is contained in a threefold of $\Delta$-genus 1,
intersection of quadrics containing it.

In this case (see \cite{Ko1}) the canonical image of $F$
is a complete intersection of two quadrics and a cubic. We follow the
notations of {\sl Case 2}. Denote $H_i=H_{\vert W_i}$. Now $\widetilde
W=W_4$ is fibred over $B$ by threefolds of degree four in ${\Proj}^5$,
complete intersections of two quadrics, and $W_2 \longrightarrow B$
is an elliptic surface over $B$.

Then we have

\bigskip

\noindent {\bf Lemma 3.3.} \
{\it

{\rm (a)} $\Phi_{\ast} {\mathcal {O}}_{\widetilde W}(H)=\pi_{\ast}
\omega_{T/B}$.

{\rm (b)} $\mbox{\rm deg}\Phi_{\ast} {\mathcal {O}}_{\widetilde W}(2H) \geq
\frac{1}{2} H^4 + 5 \mbox{\rm deg} \Phi_{\ast} {\mathcal {O}}_{\widetilde W}
(H)$.

{\rm (c)} $K^3_{T/B} \geq H^4 + 4 (\chi {\mathcal {O}}_B \chi {\mathcal {O}}_F
- \chi {\mathcal {O}}_T)$.
}

\medskip

\noindent {\petita Proof:}

(a) Follows as in Case 2.

(b) Note that, as in Case 2, formula (b) is invariant under changing
$H$ by $H + \Phi^{\ast}(A)$ so we can construct a smooth ladder of
$(\buildrel \sim \over W,H)$.
For $i \geq 2$ and $t \in B$ general $(W_i)_t \subseteq {\Proj}^{i+1}
$ is a complete intersection so it is projectively normal. On the other
side $R^1 \Phi_{\ast} {\mathcal {O}}_{W_i}$ is locally free for $i \geq 1$ (see
\cite{K2}) and in fact $R^1 \Phi_{\ast} {\mathcal {O}}_{W_i}=0$
except for $i=2$, for which it is a line bundle of degree $-\chi
{\mathcal {O}}_{W_2}$. Let $E=(W_2)_t$ any fibre of $\Phi: W_2 \longrightarrow
B$. Since $H_2+E$ is nef and big on $W_2$ we have from
Kawamata-Viehweg vanishing and the exact sequence
$$0 \rightarrow H^0(W_2,-E-H_2) \rightarrow H^0(W_2,-H_2)
\rightarrow H^0(E, {\mathcal {O}}_E(-1)) \rightarrow H^1(W_2,-E-H_2)$$
that $h^1(E,{\mathcal {O}}_E(1))=h^0(E,{\mathcal {O}}_E(-1))=0$ (recall that
$K_E={\mathcal {O}}_E$ since $W_2$ is elliptic) and hence that
$R^1\Phi_{\ast}{\mathcal {O}}_{W_2}(H)=0$. Then again by induction we have
$$\begin{array}{rl}
\forall i \geq 0 \quad \forall n \geq 1 \quad &
R^1\Phi_{\ast}{\mathcal {O}}_{W_i}(nH)=0\\
\\
\forall i \not= 2 \quad &R^1\Phi_{\ast}{\mathcal {O}}_{W_i}=0
\end{array}$$

Therefore we have exact sequences
$$\begin{array}{cl}
&0 \longrightarrow \Phi_{\ast} {\mathcal {O}}_{W_{i+1}}(H) \longrightarrow
\Phi_{\ast} {\mathcal {O}}_{W_{i+1}}(2H) \longrightarrow \Phi_{\ast}
{\mathcal {O}}_{W_i}(2H) \longrightarrow 0\quad \mbox{\rm for}\ i \geq 0\\
\\
&0 \longrightarrow \Phi_{\ast} {\mathcal {O}}_{W_{i+1}} \longrightarrow
\Phi_{\ast} {\mathcal {O}}_{W_{i+1}}(H) \longrightarrow \Phi_{\ast}
{\mathcal {O}}_{W_i}(H) \longrightarrow 0\quad \mbox{\rm for}\ i \not= 1\\
\\
&0 \longrightarrow \Phi_{\ast} {\mathcal {O}}_{W_2} \longrightarrow
\Phi_{\ast} {\mathcal {O}}_{W_2}(H) \longrightarrow \Phi_{\ast}
{\mathcal {O}}_{W_1}(H) \longrightarrow R^1 \Phi_{\ast} {\mathcal {O}}_{W_2}
\longrightarrow 0  \end{array}$$

Denote $d=\mbox{\rm deg} \Phi_{\ast} {\mathcal {O}}_{\widetilde W}(H)=
\mbox{\rm deg}\pi_{\ast} \omega_{T/C}$. Then we have
$$\begin{array}{rl}
d&=\mbox{\rm deg} \Phi_{\ast} {\mathcal {O}}_{W_4}(H)=\mbox {\rm deg}
\Phi_{\ast} {\mathcal {O}}_{W_3}(H)=\mbox{\rm deg} \Phi_{\ast}
{\mathcal {O}}_{W_2}(H)=\\
&=\mbox{\rm deg} \Phi_{\ast} {\mathcal {O}}_{W_1}(H) -\mbox{\rm deg}
R^1 \Phi_{\ast} {\mathcal {O}}_{W_2}\end{array}$$
and then
$$\mbox{\rm deg} \Phi_{\ast} {\mathcal {O}}_{\widetilde W}(2H)=4d+H^4+
\mbox{\rm deg}R^1\Phi_{\ast}{\mathcal {O}}_{W_2}=4d+H^4-\chi {\mathcal {O}}_{W_2}.$$

Note that, since $\Phi:W_2 \longrightarrow B$ is an elliptic fibration,
we have that $K_{W_2}\cong \Phi^{\ast}(L)+M$, where $M \geq 0$ and
contained in fibres and deg$L=\chi {\mathcal {O}}_{W_2}+2(b-1)$. So,
Riemann-Roch on $W_2$ and Leray spectral sequence yields
$$\begin{array}{rl}
d-4(b-1)&=\chi \Phi_{\ast} {\mathcal {O}}_{W_2}(H)=\chi {\mathcal {O}}_{W_2}(H)=
\chi {\mathcal {O}}_{W_2}+\frac{1}{2} H_2^2-\frac{1}{2} H_2 K_{W_2} \leq
\\ \\
&\leq - \chi {\mathcal {O}}_{W_2} -4 (b-1) + \frac{1}{2} H^4 \end{array}$$
since $H_2^2=H^4$, $H_2$ is nef and $H \Phi^{-1}(t)=4$ for $t \in B$.
Then $-\chi {\mathcal {O}}_{W_2} \geq d - \frac{1}{2}H^4$ and hence
deg$\Phi_{\ast} {\mathcal {O}}_{\widetilde W}(2H) \geq 5d + \frac{1}{2} H^4$.

(c) The same argument as in {\sl Case 2} works.
$\quad \Box$

\medskip

Now we only have to use $H^4\geq 0$ and get $K_{T/B}^3 \geq
4(\chi {\mathcal {O}}_B \chi {\mathcal {O}}_F- \chi {\mathcal {O}}_T)$ as
needed. Note that using good lower bounds for $H^4$ as in {\sl Case 2}
we can obtain stronger bounds for $K_{T/B}^3$ in this case.
$\quad \Box$

\bigskip

\noindent {\petita Remark 3.4.} The bounds obtained in Theorem 3.1 for
fibred canonical threefolds hold when simply $\vert K_F \vert$ induces
a birational map (it is not necessary that $T$ be canonical).

\vglue.7truecm

\noindent {\capi Appendix. The relative hyperquadrics method for
threefolds}

\vglue.3truecm

The method of counting relative hyperquadrics, originated in \cite{R2}
and \cite{CaCi} was successfully applied by Konno in
\cite{Ko4} to study the slope of fibred surfaces with small fibre
genus. Here we construct the fundamental sequence and prove the first
elementary conclusions which are needed in the previous Section.

Let $T$ be a normal, $\Rational$-factorial, projective threefold with only
terminal singularities, and let $\pi: T \longrightarrow B$ be a relatively
minimal fibration onto a smooth curve of genus $b$. Following Ohno
(\cite{O}), if $D$ is a Weil divisor on $T$ and ${\cal E}=
\pi_{\ast} {\cal O}_T(D)$ we have
$$\xymatrix{
{\widetilde T} \ar[d]_\mu \ar[dr]^{\lambda}\\
T \ar@{-->}[r]_\psi \ar[d]_\pi & Y \ar[dl]^\varphi \ar@{^{(}->}[r]^-i &
{\Proj}_B({\mathcal {E}})=:Z \\
B
}$$
where (1) $\psi$ is induced by $\pi^{\ast}\pi_{\ast} {\mathcal {O}}_T(D)
\longrightarrow {\mathcal {O}}_T(D)$ and $Y=\overline{\mbox{\rm Im}}\psi$.

(2) $\mu:\widetilde T \longrightarrow T$ is a desingularization of $T$
such that $\lambda=\psi \circ \mu$ is everywhere defined.

(3) $(\lambda^{\ast} \circ i^{\ast})L_{\mathcal {E}} \sim_{\Rational}
\mu^{\ast}(D-D_1)-E$, being $L_{\mathcal {E}}$ the tautological divisor on
$Z$, $D_1$ the codimension one base Weil divisor of ${\mathcal {O}}_T(D)$ and $E$
is an effective $\Rational$-divisor $\mu$-exceptional.

\bigskip

\noindent {\bf Proposition A-1.} Under the above hypothesis we have an
exact sequence
$$0 \longrightarrow \varphi_{\ast} {\mathcal {J}}_{Y,Z}(2 L_{\mathcal {E}})
\longrightarrow S^2\pi_{\ast} {\mathcal {O}}_T(D) \longrightarrow \pi_{\ast}
{\mathcal {O}}_T (D)^{\cite{Be1}} $$
(the generalized Max-Noether sequence associated to $\pi$).

\medskip

\noindent {\petita Proof:}

From the exact sequence
$$0 \longrightarrow
{\mathcal {J}}_{Y,Z}(2 L_{\mathcal {E}}) \longrightarrow {\mathcal {O}}_Z(2
L_{\mathcal {E}}) \longrightarrow i_{\ast} {\mathcal {O}}_Y \otimes
{\mathcal {O}}_Z(2 L_{\mathcal {E}}) \longrightarrow 0$$
we have
$$0 \longrightarrow
\varphi_{\ast} {\mathcal {J}}_{Y,Z}(2 L_{\mathcal {E}}) \longrightarrow
S^2 \pi_{\ast} {\mathcal {O}}_T(D) \longrightarrow  \varphi_{\ast}
(i_{\ast} {\mathcal {O}}_Y \otimes {\mathcal {O}}_Z(2 L_{\mathcal {E}}))$$

Now the natural map $\pi^{\ast} \pi_{\ast} {\mathcal {O}}_T(D)
\longrightarrow {\mathcal {O}}_T(D)$ induces a map
$\pi^{\ast}(S^2 \pi_{\ast} {\mathcal {O}}_T(D)) =
S^2 \pi^{\ast} \pi_{\ast} {\mathcal {O}}_T(D) \rightarrow
({\mathcal {O}}_T(D) \otimes {\mathcal {O}}_T(D))^{\ast \ast} =
{\mathcal {O}}_T(D)^{\cite{Be1}}$ and hence
$\delta : S^2 \pi _{\ast} {\mathcal {O}}_T(D) \rightarrow
\pi _{\ast} {\mathcal {O}}_T(D)^{\cite{Be1}}$. Let $K=\mbox{\rm ker}\delta$. For
general $t \in B, K_t=(\varphi_{\ast} {\mathcal {J}}_{Y,Z}(2L_
{\mathcal {E}}))_t$, and hence $K= \varphi _{\ast} {\mathcal {J}}_
{Y,Z}(2L_{\mathcal {E}})$ since $\pi_{\ast} {\mathcal {O}}_T(D)^{\cite{Be1}}$
is locally free.
$\quad \Box$

\bigskip

\noindent {\bf Corollary A-2.} Under the same hypothesis we have
$$K_{T/B}^3 \geq (2 p_g(F)-4)(\chi {\mathcal {O}}_B \chi {\mathcal {O}}_F
-\chi {\mathcal {O}}_T)-2 \mbox{\rm deg}K -2 \ell (2) \qquad \qquad
\qquad \ (1)$$
where $K=\varphi_{\ast} {\mathcal {J}}_{Y,Z}(2 L_{\mathcal {E}})$
$F=\pi^{-1}(t)$ for $t \in B$, and {\it l}(2) is the second order
correction term of Reid-Fletcher to the plurigenera of $T$ (cf. \cite{F}).

\medskip

\noindent {\petita Proof:}

Let $D=K_{T/B}$ (which is in general a Weil divisor) and take
degrees in the generalized Max-Noether sequence. Use
$$d=\mbox {\rm deg}\pi_{\ast} {\omega}_{T/B} \geq
(\chi {\mathcal {O}}_B \chi {\mathcal {O}}_F -
\chi {\mathcal {O}}_T) \qquad
\qquad \qquad \mbox{\rm (\cite{O} p. 656})$$
$$\mbox{\rm deg}\pi_{\ast} \omega_{T/B}^{\cite{Be1}}=\frac{1}{2}
K_{T/B}^3+3(\chi {\mathcal {O}}_B \chi {\mathcal {O}}_F -\chi {\mathcal {O}}_T)
+ \ell (2) \qquad \mbox{\rm (\cite{O} Lemma 2.8)}$$
$$\mbox{\rm deg}S^2\pi_{\ast} \omega_{T/B}=(p_g(F) +1)d$$
$$\mbox{\rm rk}S^2\pi_{\ast} \omega_{T/B} =
{{p_g(F)+1}\choose{2}}$$
and that if ${\mathcal {C}}=$coker($S^2 \pi _{\ast} \omega_{T/B}
\rightarrow \pi_{\ast} \omega_{T/B}^{\cite{Be1}}$),
deg${\mathcal {C}} \geq 0$ since $\pi_{\ast} \omega_{T/B}^{\cite{Be1}}$ is
semipositive (\cite{O}).
$\quad \Box$

\bigskip

\noindent {\petita Remark A.3.} For small values of the invariants
$p_g(F)$, $q(F)$, $K_F^2$ it could be interesting to consider
$D=m K_{T/B}$ for $m \geq 1$. We obtain then bounds for
$K_{T/B}^3$ which are better than (1).

\bigskip

In general deg$K$ is difficult to be computed or bounded. There are some
special cases where this is easier. Notice that rk$K=h^0(I_{\Sigma,
{\proj}^r}(2))$ where $\Sigma$ is the canonical image of $F$ and
$r=p_g(F)-1$. Then following Lemma 1.1 we have
$$\begin{array}{ll}
h^0({\mathcal {J}}_{\Sigma, {\proj}^r}(2)) \leq \frac{(r-2)(r-3)}{2}
&\qquad \mbox{\rm if}\ \Sigma \ \mbox{\rm is a non ruled surface}\\
\\
h^0({\mathcal {J}}_{\Sigma, {\proj}^r}(2)) \leq \frac{(r-1)(r-2)}{2}-q(\Sigma)
&\qquad \mbox{\rm if}\ \Sigma \ \mbox{\rm is a ruled surface}\\
\\
h^0({\mathcal {J}}_{\Sigma, {\proj}^r}(2)) \leq \frac{r(r-1)}{2}
&\qquad \mbox{\rm if}\ \Sigma \ \mbox{\rm is a curve}\end{array}$$

\bigskip

\noindent {\bf Lemma A.4.}
{\it

{\rm (a)} If $p_g(F) \geq 2$ and ${\mathcal {E}}=\pi_{\ast} \omega_{T/B}$ is semistable
then {\rm deg}$K \leq 2 \frac{\mbox{\rm rk} \mbox{\it K}}{p_g(F)}d$.

{\rm (b)} If $K={\mathcal {L}}_1 \oplus \dots {\mathcal {L}}_s$ (s=\mbox {\rm rk}K)
then
$$\mbox{\rm deg}K \leq (\mbox{{\rm rk}{\it K}})\frac{2}{3}d$$
(in particular this happens if $s \leq 1$ or $b=0$).
}

\medskip

\noindent {\petita Proof:}

(a) If $\mathcal {E}$ is semistable then so it is $S^2\mathcal {E}$. Then
we use the natural inclusion $K \hookrightarrow S^2\mathcal {E}$.

(b) If $x_i=\mbox{\rm deg}{\mathcal {L}}_i$ then there exists a section
$s \in H^0(K \otimes {\mathcal {L}}_i^{-1})\cong
H^0({\mathcal {J}}_{Y,Z}(2 L_{\mathcal {E}})\otimes
{\mathcal {O}}_Z(\varphi^{\ast}({\mathcal {L}}_i^{-1}))
\hookrightarrow H^0(Z, {\mathcal {O}}(2 L_{\mathcal {E}})\otimes
{\varphi}^{\ast}({\mathcal {L}}_i^{-1}))$ so there exists a relative
hyperquadric $Q_i \equiv 2 L_{\mathcal {E}}-x_i {\varphi}^{-1}(t)$
(numerical equivalence). The result follows then from the following Lemma
which is a slight refinement of \cite{Ko4} Remark 1.7, and the
fact that for every $i$, rk${Q}_i \geq 3$.
$\quad \Box$

\medskip

\noindent {\bf Lemma A-5.}
{\it
Let $Q \equiv 2 L_{\mathcal {E}}-x{\varphi}^{-1}(t)$ be a relative
hyperquadric. Let $\nu_1 \geq \nu_2 \geq \dots \geq \nu_k$ the virtual
slopes of the Harder-Narashiman filtration of $\mathcal {E}$
($k=\mbox{\rm rk}\mathcal {E}$). Let $p=\mbox{\rm rk} Q$; then
$$x \leq \mbox{\rm min}_{1 \leq i \leq p} \{\nu_i+\nu_{p-i}\} \leq
\frac{2}{p} \mbox{\rm deg}\mathcal {E} $$
}

\bigskip

\noindent {\bf Corollary A.6.}
{\it
With the same notations as above, assume $p_g(F)\geq 2$.

{\rm (a)} If ${\mathcal {E}} =\pi_{\ast} \omega_{T/B}$ is
semistable then
$$\begin{array}{rl}
K_{T/B}^3 \geq \left (10 -\frac{24}{p_g(F)}\right )
(\chi {\mathcal {O}}_B \chi {\mathcal {O}}_F - \chi {\mathcal {O}}_T) -2 \ell (2) &\quad
\mbox {if}\ \Sigma \ \mbox{is a non-ruled surface}\\
K_{T/B}^3 \geq \left (6 -\frac{12}{p_g(F)}\right )
(\chi {\mathcal {O}}_B \chi {\mathcal {O}}_F - \chi {\mathcal {O}}_T) -2 \ell (2) &\quad
\mbox {if}\ \Sigma \ \mbox{is a ruled surface}\\
K_{T/B}^3 \geq \left (2 -\frac{4}{p_g(F)}\right )
(\chi {\mathcal {O}}_B \chi {\mathcal {O}}_F - \chi {\mathcal {O}}_T) -2 \ell (2) &\quad
\mbox {if}\ \Sigma \ \mbox{is a curve}
\end{array}$$

{\rm (b)} If $h^0({\mathcal {J}}_{\Sigma,{\proj}^r}(2))=0$ then
$$K^3_{T/B} \geq  (2 p_g(F)-4)(\chi {\mathcal {O}}_F \chi {\mathcal {O}}_B
- \chi {\mathcal {O}}_T) - 2 \ell(2)$$

{\rm (c)} If $h^0({\mathcal {J}}_{\Sigma,{\proj}^r})(2)=1$ then
$$K^3_{T/B} \geq  (2 p_g(F)-\frac{16}{3})(\chi {\mathcal {O}}_F
\chi {\mathcal {O}}_B - \chi {\mathcal {O}}_T) - 2 \ell(2)$$
}

\newpage

\noindent {\petita Proof:}

Take degrees at the generalized Max-Noether sequence and use Remark
A.3 and Lemma A.4.
$\quad \Box$

\vglue.7truecm

\vglue.5truecm

Departament de Matematica Aplicada I

ETSEIB. UNIVERSITAT POLITECNICA DE CATALUNYA

Avda. Diagonal 647

08190-BARCELONA

SPAIN

e-mail: barja@ma1.upc.es

\end{document}